\documentclass{article}
\usepackage{amsmath,amssymb,amscd,enumerate,bbm,geometry}

\catcode`\à=\active \defà{\`a} \catcode`\À=\active \defÀ{\`A}
\catcode`\á=\active \defá{\'a} \catcode`\Á=\active \defÁ{\'A}
\catcode`\ä=\active \defä{\"a} \catcode`\Ä=\active \defÄ{\"A}
\catcode`\â=\active \defâ{\^a} \catcode`\Â=\active \defÂ{\^A}
\catcode`\å=\active \defå{{\aa}} \catcode`\Å=\active \defÅ{{\AA}}
\catcode`\ç=\active \defç{\c{c}} \catcode`\Ç=\active \defÇ{\c{C}}
\catcode`\è=\active \defè{\`e} \catcode`\È=\active \defÈ{\`E}
\catcode`\é=\active \defé{\'e} \catcode`\É=\active \defÉ{\'E}
\catcode`\ë=\active \defë{\"e} \catcode`\Ë=\active \defË{\"E}
\catcode`\ê=\active \defê{\^e} \catcode`\Ê=\active \defÊ{\^E}
\catcode`\ì=\active \defì{\`{\i}} \catcode`\Ì=\active \defÌ{\`{\I}}
\catcode`\í=\active \defí{\'{\i}} \catcode`\Í=\active \defÍ{\'{\I}}
\catcode`\ï=\active \defï{\"{\i}} \catcode`\Ï=\active \defÏ{\"{\I}}
\catcode`\î=\active \defî{\^{\i}} \catcode`\Î=\active \defÎ{\^{\I}}
\catcode`\ò=\active \defò{\`o} \catcode`\Ò=\active \defÒ{\`O}
\catcode`\ó=\active \defó{\'o} \catcode`\Ó=\active \defÓ{\'O}
\catcode`\ö=\active \defö{\"o} \catcode`\Ö=\active \defÖ{\"O}
\catcode`\ô=\active \defô{\^o} \catcode`\Ô=\active \defÔ{\^O}
\catcode`\ù=\active \defù{\`u} \catcode`\Ù=\active \defÙ{\`U}
\catcode`\ú=\active \defú{\'u} \catcode`\Ú=\active \defÚ{\'U}
\catcode`\ü=\active \defü{\"u} \catcode`\Ü=\active \defÜ{\"U}
\catcode`\û=\active \defû{\^u} \catcode`\Û=\active \defÛ{\^U}
\catcode`\ý=\active \defý{\'y} \catcode`\Ý=\active \defÝ{\'Y}
\catcode`\ÿ=\active \defÿ{\"y} \catcode`\˜=\active \def˜{\"Y}
\catcode`\½=\active \def½{!`}
\catcode`\¾=\active \def¾{?`}
\catcode`\ß=\active \defß{{\ss}}
\newcommand{\email}[1]{{{\em Email}: {\tt #1}}}
\newcommand{\tmop}[1]{\ensuremath{\operatorname{#1}}}
\newenvironment{enumeratenumeric}{\begin{enumerate}[1.]}{\end{enumerate}}
\newtheorem{definition}{Definition}
\newenvironment{enumeratealpha}{\begin{enumerate}[a{\textup{)}}]}{\end{enumerate}}
\newtheorem{proposition}{Proposition}
\newenvironment{proof}{
  \noindent\textbf{Proof}\ }{\hspace*{\fill}
  \begin{math}\Box\end{math}\medskip}
\newtheorem{notation}{Notation}
\newenvironment{itemizeminus}
  {\begin{itemize}}{\end{itemize}}
\newcommand{\tmmathbf}[1]{\ensuremath{\boldsymbol{#1}}}
\newtheorem{theorem}{Theorem}

\begin{document}

\title{Formulas for Birkhoff-(Rota-Baxter) decompositions related to connected
bialgebras.} 
\author{Frédéric Menous\footnote{Univ Paris-Sud, Laboratoire de Math\'ematiques d'Orsay, Orsay
  Cedex, F-91405; CNRS, Orsay cedex,
  F-91405. \email{Frederic.Menous@math.u-psud.fr}}}
\date{}

\maketitle

\begin{abstract}
  In recent years, The usual BPHZ algorithm for renormalization in quantum
  field theory has been interpreted, after dimensional regularization, as the
  Birkhoff-(Rota-Baxter) decomposition (BRB) of characters on the Hopf algebra
  of Feynman graphs, with values in a Rota-Baxter algebra.
  
  We give in this paper formulas for the BRB decomposition in the group
  $\mathcal{C}( H, A )$ of characters on a connected Hopf algebra $H$, with
  values in a Rota-Baxter (commutative) algebra $A$.
  
  To do so we first define the stuffle (or quasi-shuffle) Hopf algebra
  $A^{\tmop{st}}$ associated to an algebra $A$. We prove then that for any
  connected Hopf algebra $H = k 1_H \oplus H'$, there exists a canonical
  injective morphism from $H$ to $H'^{\tmop{st}}$. This morphism induces an
  action of $\mathcal{C}( A^{\tmop{st}}, A )$ on $\mathcal{C}( H, A )$ so that
  the BRB decomposition in $\mathcal{C}( H, A )$ is determined by the action
  of a unique (universal) element of $\mathcal{C}( A^{\tmop{st}}, A )$.
\end{abstract}

\section{Introduction.}

In this paper we deal with connected bialgebras $H = k 1_H \oplus H'$. As
reminded in section \ref{sec2}, such bialgebras are Hopf algebras and the
coalgebra structure on $H$ induces a convolution product on the space
$\mathcal{L}( H, A )$ of linear morphisms from $H$ to an associative algebra
$A$. If $A$ is unital, then the subset $\mathcal{U}( H, A )$ of linear
morphisms that send $1_H$ on $1_A$ is a group for the convolution and, if A is
commutative the subset $\mathcal{C}( H, A )$ of characters (algebra morphisms)
is a subgroup of $\mathcal{U}( H, A )$.

In section \ref{sec3}, the target  unital algebra $A$ is equipped with a
Rota-Baxter operator $p_+$ :
\begin{enumeratenumeric}
  \item $p_+^2 = p_+$.
  
  \item $A = \tmop{Im} p_+ \oplus \tmop{Im} p_- = A_+ \oplus A_-$ ($p_- =
  \tmop{Id} - p_+$).
  
  \item $A_+$ and $A_-$ are subalgebras.
\end{enumeratenumeric}
With this hypothesis, It is well-know that there exists a unique
Birkhoff-(Rota-Baxter) (or BRB) decomposition of any morphism $\varphi \in
\mathcal{U}( H, A )$
\[ \varphi_- \ast \varphi = \varphi_+ \hspace{2em} \varphi_+, \varphi_- \in
   \mathcal{U}( H, A ) \]
where $\varphi_+ ( H' ) \subset A_+$ and $\varphi_+ ( H' ) \subset A_-$.
Moreover, if $A$ is commutative, this decomposition is defined in the subgroup
$\mathcal{C}( H, A )$. The proof of this result is recursive, using the
filtration on $H$. We propose to give explicit, and in some sense universal,
formulas for $\varphi_+$ and $\varphi_-$.

To do so, we define in section \ref{sec4} the stuffle Hopf algebra
$A^{\tmop{st}}$ of an algebra $A$, that is, the standard tensor coalgebra of
$A$ equipped with the stuffle (or quasi-shuffle product). With give then the
example of computation for the inverse and the BRB decomposition of a
fundamental element $j \in \mathcal{U}( A^{\tmop{st}}, A )$ defined by
\[ j ( \emptyset ) = 1_A, \hspace{1em} j ( a_1 ) = a_1, \hspace{1em} j ( a_1
   \otimes \ldots \otimes a_s ) = 0 \hspace{1em} \tmop{if} \hspace{1em} s \geq
   2 \]

Stuffle Hopf algebras play a central role since, as proved in section
\ref{sec5}, For any given connected Hopf algebra $H = k 1_H \oplus H'$ there
exists a canonical Hopf morphism $\iota : H \rightarrow H'^{\tmop{st}}$
defined with the help of the reduced coproduct. We prove then in section
\ref{sec6} that this morphism induces an action of $\mathcal{U}(
A^{\tmop{st}}, A )$ on $\mathcal{U}( H, A )$. More precisely we define a map
$T :\mathcal{U}( A^{\tmop{st}}, A ) \times \mathcal{U}( H, A ) \rightarrow
\mathcal{U}( H, A )$ such that
\[ T ( j, \varphi ) = \varphi \hspace{1em} \tmop{and} \hspace{1em} T ( f \ast
   g, \varphi ) = T ( f, \varphi ) \ast T ( g, \varphi ) \]
and we obtain explicit formulas as follows :
\begin{enumeratenumeric}
  \item If $j^{\ast - 1}$ is the inverse of $j$, then $\varphi^{\ast - 1} = T
  ( j^{\ast - 1}, \varphi )$.
  
  \item If $j_- \ast j = j_+$ (BRB decomposition), the $\varphi_- \ast \varphi
  = \varphi_+$ where $\varphi_{\pm} = T ( j_{\pm}, \varphi$).
\end{enumeratenumeric}

Note that such formulas were already derived in the Hopf algebra of Feynman
graphs using matrix calculus (see {\cite{mat}}) but the formulas presented
here are intrinsic : no choice of a basis is needed since we don't have to deal
with matrices.

\section{Reminder about connected bialgebras.} \label{sec2}

We follow here the notations and definitions given in {\cite{fig}} (see also
{\cite{maj}} and {\cite{sweed}}). Let $k$ be a commutative field.

\subsection{Bialgebras and Hopf algebras.}

In the sequel, we will work with  bialgebras over $k$:

\begin{definition}
  A bialgebra H is a $k$--vector space equipped with four linear maps $m : H
  \otimes H \rightarrow H$ (product : $m ( x \otimes y ) = x y$) , $u : k
  \rightarrow H$ (unit : $u ( 1_k ) = 1_H$), $\Delta : H \rightarrow H \otimes
  H$ (coproduct) and $\eta : H \rightarrow k$ (counit) such that
  \begin{enumeratenumeric}
    \item $( H, m, u )$ is a unital associative algebra. This reads :
    \begin{enumeratealpha}
      \item Associativity : $m \circ ( m \otimes \tmop{Id} ) = m \circ (
      \tmop{Id} \otimes m ) : H \otimes H \rightarrow H$;
      
      \item Unit :  $m \circ ( u \otimes \tmop{Id} ) = m \circ ( \tmop{Id}
      \otimes u ) = \tmop{Id} : k \otimes H = H \otimes k = H \rightarrow H$.
    \end{enumeratealpha}
    \item $( H, \Delta, \eta )$ is a coassociative coalgebra with a counit :
    \begin{enumeratealpha}
      \item Coassociativity : $( \Delta \otimes \tmop{Id} ) \circ \Delta = (
      \tmop{Id} \otimes \Delta ) \circ \Delta : H \rightarrow H \otimes H
      \otimes H$ ;
      
      \item Counit : $( \eta \otimes \tmop{Id} ) \circ \Delta = ( \tmop{Id}
      \otimes \eta ) \circ \Delta = \tmop{Id} : H \rightarrow H$.
    \end{enumeratealpha}
    \item The following diagram commutes :
    \[ \begin{CD}
 H \otimes H @>m>> H \\
@VV{\Delta \otimes \Delta}V @AA{m \otimes m}A \\
H \otimes H \otimes H \otimes H @>{\tmop{Id} \otimes \tau
         \otimes \tmop{Id}}>> H \otimes H \otimes H \otimes H
\end{CD}
\]


    where $\tau$ is the linear map defined by $\tau ( h \otimes g ) = g
    \otimes h$ and
    \[ \Delta ( 1_H ) = 1_H \otimes 1_H, \hspace{1em} \eta ( h g ) = \eta ( h
       ) \eta ( g ) \]
  \end{enumeratenumeric}
\end{definition}

Note that $p = u \circ \eta$ is an idempotent ($p' = \tmop{Id} - p )$ and $H =
\tmop{Im} \hspace{1em} p \oplus \tmop{Ker} \hspace{1em} p = \tmop{Im}
\hspace{1em} u \oplus \tmop{Ker} \hspace{1em} \eta = k 1_H \oplus H'$.

If $\Delta h = \sum_j h_{j ( 1 )} \otimes h_{j ( 2 )}$, we will sometimes
write
\begin{equation}
  \Delta h = \sum h_{( 1 )} \otimes h_{( 2 )} = h_{( 1 )} \otimes h_{( 2 )}
\end{equation}
For example, coassociativity read
\[ \sum h_{( 1 ) ( 1 )} \otimes h_{( 1 ) ( 2 )} \otimes h_{( 2 )} = \sum h_{(
   1 )} \otimes h_{( 2 ) ( 1 )} \otimes h_{( 2 ) ( 2 )} = \sum h_{( 1 )}
   \otimes h_{( 2 )} \otimes h_{( 3 )} \]
and thanks to coassociativity, we can define recursively and without any
ambiguity the linear morphisms $\Delta^{[ n ]} : H \rightarrow H^{\otimes^n}$
($n \geq 1$) by $\Delta^{[ 1 ]} = \tmop{Id}$ and, for $n \geq 1$,
\begin{equation}
  \Delta^{[ n + 1 ]} = ( \tmop{Id} \otimes \Delta^{[ n ]} ) \circ \Delta = (
  \Delta^{[ n ]} \otimes \tmop{Id} ) \circ \Delta = ( \Delta^{[ k ]} \otimes
  \Delta^{[ n + 1 - k ]} ) \circ \Delta \hspace{1em} ( 1 \leq k \leq n )
\end{equation}
and write
\begin{equation}
  \Delta^{[ n ]} h = \sum h_{( 1 )} \otimes \ldots \otimes h_{( n )}
\end{equation}
On the same way, for $n \geq 1$, we define $m^{[ n ]} : H^{\otimes^n}
\rightarrow H$ by $m^{[ 1 ]} = \tmop{Id}$ and
\begin{equation}
  m^{[ n + 1 ]} = m \circ ( \tmop{Id} \otimes m^{[ n ]} ) = m \circ ( m^{[ n
  ]} \otimes \tmop{Id} )
\end{equation}
Such a bialgebra is a Hopf algebra if there exists an antipode $S$, that is to
say a linear map $S : H \rightarrow H$ such that :
\begin{equation}
  m \circ ( \tmop{Id} \otimes S ) \circ \Delta = m \circ ( S \otimes \tmop{Id}
  ) \circ \Delta = u \circ \eta : H \rightarrow H
\end{equation}
From now on, we should focus on connected bialgebras since they are
automatically Hopf algebras.

\subsection{Connected bialgebras.}

\begin{definition}
  A bialgebra $\mathbbm{N}$--filtered as a vector space is called a filtered
  bialgebra when the filtering is compatible with both the algebra an the
  coalgebra structure; that is, there exist a sequence of subspaces $H_0
  \varsubsetneq H_1 \varsubsetneq \ldots$ such that $\bigcup_{n \geq 0} H_n =
  H$ and
  \begin{equation}
    \Delta H_n \subseteq \bigoplus_{k = 0}^n H_k \otimes H_{n - k} ;
    \hspace{1em} H_n H_m \subseteq H_{n + m}
  \end{equation}
  Connected bialgebras are those filtered bialgebras such that $H_0 = k 1_H =
  \tmop{Im} u = \tmop{Im} p$.
\end{definition}

We shall come back in the next section on the fact that the antipode
automatically exists. For such a connected bialgebra, if, for all $n \geq 1$,
$H'_n = H' \cap H_n$, then
\begin{equation}
  \forall h \in H'_n, \hspace{1em} \Delta h = 1 \otimes h + h \otimes 1 + y
  \hspace{1em} \tmop{where} \hspace{1em} y = ( p' \otimes p' ) \circ \Delta h
  \in \bigoplus_{k = 1}^{n - 1} H'_k \otimes H'_{n - k}
\end{equation}
and, it will be useful to define the reduced coproduct $\Delta'$ on $H'$
defined by
\begin{equation}
  \Delta' h = \Delta h - 1 \otimes h - h \otimes 1
\end{equation}
Then, recursively, for $n \geq 1$, we define $\Delta'^{[ n ]} = p'^{\otimes^n}
\circ \Delta^{[ n ]} : H' \rightarrow H'^{\otimes^n}$. We note
\begin{equation}
  \forall h \in H', \hspace{1em} \Delta'^{[ n ]} h = \sum h'_{( 1 )} \otimes
  \ldots \otimes h'_{( n )}
\end{equation}
and we have, for $h$ in $H$, $p'^{\otimes^n} \circ \Delta^{[ n ]} ( h ) =
\Delta'^{[ n ]} \circ p' ( h )$.

Moreover, for $n \geq 1$ and $h \in H'$
\begin{equation}
  \Delta'^{[ n + 1 ]} ( h ) = ( \Delta'^{[ k ]} \otimes \Delta'^{[ n + 1 - k
  ]} ) \circ \Delta' ( h ) \hspace{1em} ( 1 \leq k \leq n )
\end{equation}
and if $h \in H'_k$ ($k \geq 1$) then, for $n > k$, $\Delta'^{[ n ]} ( h ) =
0$.

Given a connected bialgebra $H$ and an algebra $( A, m_A, u_A )$ the coalgebra
structure of $H$ induces an associative convolution product on the vector
space $\mathcal{L}( H, A )$ of $k$--linear maps :
\begin{equation}
  \forall ( f, g ) \in \mathcal{L}( H, A ) \times \mathcal{L}( H, A ),
  \hspace{1em} f \ast g = m_A \circ ( f \otimes g ) \circ \Delta
\end{equation}
with a unit given by $u_A \circ \eta$, such that $(\mathcal{L}( H, A ), \ast,
u_A \circ \eta )$ is an associative unital algebra.

\subsection{The group $(\mathcal{U}( H, A ), \ast )$.}

\begin{proposition}
  Let
  \begin{equation}
    \mathcal{U}( H, A ) = \{ f \in \mathcal{L}( H, A ) \hspace{1em} ;
    \hspace{1em} f ( 1_H ) = 1_A \}
  \end{equation}
  then $\mathcal{U}( H, A )$ is a group for the convolution product.
\end{proposition}

\begin{proof}
  $\mathcal{U}( H, A )$ is obviously stable for the convolution product and
  following {\cite{fig}} we should remind why any element $f \in \mathcal{U}(
  H, A )$ as a unique inverse $f^{\ast - 1}$ in $\mathcal{U}( H, A )$. There
  are two ways to define this inverse.
  
  Since $H = H_0 \oplus H'$ it is sufficient to define recursively $f^{\ast -
  1}$ on $H'_n$ for $n \geq 1$ ($f^{- 1} ( 1_H ) = 1_A )$. For $n \geq 1$ and
  $h \in H'_n$, we have
  \[ \Delta h = 1 \otimes h + h \otimes 1 + \sum h'_{( 1 )} \otimes h'_{( 2 )}
  \]
  where
  \[ \sum h'_{( 1 )} \otimes h'_{( 2 )} \in \bigoplus_{k = 1}^{n - 1} H'_k
     \otimes H'_{n - k} \]
  thus
  \begin{equation}
    f^{\ast - 1} \ast f ( h ) = u_A ( \eta ( a ) ) = 0 = f^{\ast - 1} ( h ) +
    f ( h ) + \sum f^{\ast - 1} ( h'_{( 1 )} ) f ( h'_{( 2 )} )
  \end{equation}
  and this defines $f^{\ast - 1}$ recursively. On the other hand, we can write
  \begin{equation}
    f^{\ast - 1} = ( u_A \circ \eta - ( u_A \circ \eta - f ) )^{\ast - 1} =
    u_A \circ \eta + \sum_{k \geq 1} ( u_A \circ \eta - f )^{\ast^k}
  \end{equation}
  In fact, this series seems to be infinite but for $h \in H'_n$
  \begin{equation}
    ( u_A \circ \eta - f )^{\ast^k} ( h ) = ( - 1 )^k \sum f ( h'_{( 1 )} )
    \ldots f ( h'_{( k )} ) = ( - 1 )^k m_A^{[ k ]} \circ f^{\otimes^k} \circ
    \Delta'^{[ k ]} ( h ) \label{eq15}
  \end{equation}
  vanishes as soon as $k > n$.

\end{proof}

The principle of recursive computation will be useful when dealing with
Birkhoff-Rota-Baxter decomposition and the main goal of this paper will be to
find also formulas like \ref{eq15}.

\begin{notation}
  If $B \subset A$ is a subalgebra of $A$ which is not unital, then we write
  \[ \mathcal{U}( H, B ) = \{ f \in \mathcal{L}( H, A ) \hspace{1em} ;
     \hspace{1em} f ( 1_H ) = 1_A \hspace{1em} \tmop{and} \hspace{1em} f ( H'
     ) \subset B \} \]
  This is a subgroup of $\mathcal{U}( H, A )$.
\end{notation}

If this result is applied to $\tmop{Id} : H \rightarrow H \in
\text{$\mathcal{U}( H, H )$}$, then its convolution inverse is the antipode
$S$ and this proves that any connected bialgebra is a Hopf algebra. Moreover
$S$ is an antialgebra morphism :
\begin{equation}
  S ( g h ) = S ( h ) S ( g )
\end{equation}
\[  \]

\subsection{Algebra morphisms or characters.}

Let $\mathcal{C}( H, A )$ the subset of $\mathcal{L}( H, A )$ whose elements
are algebra morphisms (also called characters over $A$). Of course,
\[ \mathcal{C}( H, A ) \subset \text{$\mathcal{U}( H, A )$} \]
but this shall not be a subgroup. If $A$ is not commutative, there is no
reason why this should be stable for the convolution product. Moreover if $f
\in \text{$\mathcal{U}( H, A )$}$ is an algebra map, then its inverse $f^{\ast
- 1}$ in $\mathcal{U}( H, A )$ is an antialgebra map. In fact, if $f$ is an
algebra map then $f^{\ast - 1} = f \circ S$ :
\begin{equation}
  \begin{array}{ccc}
    f \ast f \circ S & = & m_A \circ ( f \otimes f \circ S ) \circ \Delta\\
    & = & m_A \circ ( f \otimes f ) \circ ( \tmop{Id} \otimes S ) \circ
    \Delta\\
    & = & f \circ m \circ ( \tmop{Id} \otimes S ) \circ \Delta\\
    & = & f \circ u \circ \eta\\
    & = & u_A \circ \eta
  \end{array}
\end{equation}
Nonetheless if $A$ is commutative, then $\mathcal{C}( H, A )$ is a subgroup of
$\text{$\mathcal{U}( H, A )$}$.

\section{Rota-Baxter algebras and Birkhoff-type decomposition in $\mathcal{U}(
H, A )$.} \label{sec3}

Following {\cite{kur}}, let $p_+$ an idempotent of $\mathcal{L}( A, A )$ where
$A$ is a unital algebra. If we have for $x, y$ in $\mathcal{A}$ :

\begin{equation}
  p_+ ( x ) p_+ ( y ) + p_+ ( x y ) = p_+ ( x p_+ ( y ) ) + p_+ ( p_+ ( x ) y
  ) )
\end{equation}
Then $p_+$ is a Rota-Baxter operator, $( A, p_+ )$ is a Rota-Baxter algebra
and if $p_- = \tmop{Id} - p_+$, $A_+ = \tmop{Im} p_+$ and $A_- = \tmop{Im}
p_-$ then
\begin{itemizeminus}
  \item $A = A_+ \oplus A_-$.
  
  \item $p_-$ satisfies the same relation.
  
  \item $A_+$ and $A_-$ are subalgebras.
\end{itemizeminus}
Conversely if $A = A_+ \oplus A_-$ and $A_+$ and $A_-$ are subalgebras, then
the projection $p_+$ on $A_+$ parallel to $A_-$ defines a Rota-Baxter algebra
$( A, p_+ )$.

The principle of renormalization in physics can be formulated in the following
way

\begin{proposition}
  Let $H$ be a connected bialgebra and $( A, p_+ )$ a Rota-Baxter algebra then
  for any $\varphi \in \mathcal{U}( H, A )$ there exists a unique pair $(
  \varphi_+, \varphi_- ) \in \mathcal{U}( H, A_+ ) \times \mathcal{U}( H, A_-
  )$ such that
  \begin{equation}
    \varphi_- \ast \varphi_{} = \varphi_+
  \end{equation}
  Moreover, if $A$ is commutative and $\varphi$ is a character over A, then
  $\varphi_+$ and $\varphi_-$ are also characters. This factorization will be
  called the Birkhoff-Rota-Baxter (or BRB) decomposition of $\varphi$.
\end{proposition}

\begin{proof}
  We will give the proof for characters later. As $A_+$ and $A_-$ are
  subalgebras of $A$, $\mathcal{U}( H, A_+ )$ and $\mathcal{U}( H, A_- )$ are
  subgroups of $\mathcal{U}( H, A )$.
  
  If such a factorization exists, then it is unique : If $\varphi =
  \varphi_-^{\ast - 1} \ast \varphi_+ = \psi_-^{\ast - 1} \ast \psi_+$, then
  \[ \phi = \psi_+^{} \ast \varphi^{\ast - 1}_+ = \psi_- \ast \varphi_-^{\ast
     - 1} \in \mathcal{U}( H, A_+ ) \cap \mathcal{U}( H, A_- ) \]
  thus for $h \in H'$, $\phi ( h ) \in A_+ \cap A_- = 0$. We finally get that
  \[ \psi_+^{} \ast \varphi^{\ast - 1}_+ = \psi_- \ast \varphi_-^{\ast - 1} =
     u_A \circ \eta \]
  and $\varphi_+ = \psi_+$, $\varphi_- = \psi_-$.
  
  Let us prove now that such a factorization exists. Let $\varphi \in
  \mathcal{U}( H, A )$, we must have $\varphi_+ ( 1_H ) = \varphi_- ( 1_H ) =
  1_A$. Let $\bar{\varphi} \in \mathcal{U}( H, A )$ the Bogoliubov preparation
  map defined recursively on vector spaces $H'_n$ ($n \geq 1$) by
  \begin{equation}
    \bar{\varphi} ( h ) = \varphi ( h ) - m_A \circ ( p_- \otimes \tmop{Id}
    ) \circ ( \bar{\varphi} \otimes \varphi ) \circ \Delta' ( h )
  \end{equation}
  Now if $\varphi_+$ and $\varphi_-$ are the elements of $\mathcal{U}( H, A )$
  defined on $H'$ by
  \[ \varphi_+ ( h ) = p_+ \circ \bar{\varphi} ( h ) \hspace{1em},
     \hspace{1em} \varphi_- ( h ) = - p_- \circ \bar{\varphi} ( h )
     \hspace{1em} ( \bar{\varphi} ( h ) = \varphi_+ ( h ) - \varphi_- ( h ) )
  \]
  Then it is clear that
  \[ \text{} \varphi_+ \in \mathcal{U}( H, A_+ ) \hspace{1em}, \hspace{1em}
     \varphi_- \in \text{$\mathcal{U}( H, A_- )$} \hspace{1em}, \hspace{1em}
     \varphi_- \ast \varphi_{} = \varphi_+ \]
  
\end{proof}

\section{The stuffle Hopf algebra $A^{\tmop{st}}$ of an algebra $A.$}
\label{sec4}

For details on the stuffle (or quasi-shuffle) product, the reader can refer to
{\cite{hof}}.

\subsection{Definition and properties.}

Let $A$ be an associative algebra. $A^{\tmop{st}}$ is the graded vector space
$A^{\tmop{st}} = \bigoplus_{n \geq 0} A^{\tmop{st}}_{( n )}$ where, for $n
\geq 1$, $A^{\tmop{st}}_{( n )} = A^{\otimes^n}$ and  $A^{\tmop{st}}_{( 0 )} =
k \emptyset$ where $\emptyset$ is a symbol for the empty tensor product. It is
obviously graded and we note $l (\tmmathbf{a}) = n$ the length of an element
$\tmmathbf{a}$ of $A^{\tmop{st}}_{( n )}$. For convenience, an element
$\tmmathbf{a}= a_1 \otimes \ldots \otimes a_r$ of $A^{\tmop{st}}$ we be called
a tuple or a word and if $\tmmathbf{a}$ and $\tmmathbf{b}$ are two words, then
$\tmmathbf{a} \otimes \tmmathbf{b}$ is the concatenation of the words. Note
also that, as $\otimes$already denotes the tensor product in $A$, when there
may be some ambiguity, we use $\otimes_{\tmop{st}}$ for the tensor product of
elements of $A^{\tmop{st}}$.

One can define recursively the stuffle or quasi shuffle product
$m_{\tmop{st}} : A^{\tmop{st}} \otimes_{\tmop{st}} A^{\tmop{st}} \rightarrow
A^{\tmop{st}}$ on $A^{\tmop{st}}$:
\begin{enumeratenumeric}
  \item For any $\tmmathbf{a} \in A^{\tmop{st}}$, $m_{\tmop{st}} ( \emptyset
  \otimes_{\tmop{st}} \tmmathbf{a}) = m_{\tmop{st}} ( \tmmathbf{a}
  \otimes_{\tmop{st}} \emptyset ) = \tmmathbf{a}$
  
  \item Let $\tmmathbf{a}= a_1 \otimes \ldots \otimes a_r \in A^{\tmop{st}}_{(
  r )}$ and $\tmmathbf{b} = b_1 \otimes \ldots \otimes b_s \in
  A^{\tmop{st}}_{( s )}$ with $r \geq 1$ and $s \geq 1$. If
  $\tilde{\tmmathbf{a}} = a_1 \otimes \ldots \otimes a_{r - 1} \in
  A^{\tmop{st}}_{( r - 1 )}$ ($\tilde{\tmmathbf{a}} = \emptyset$ if $r = 1$)
  and $\tilde{\tmmathbf{b}} = b_1 \otimes \ldots \otimes b_{s - 1} \in
  A^{\tmop{st}}_{( s - 1 )}$ ($\widetilde{\tmmathbf{b}} = \emptyset$ if $r =
  1$), then :
  \begin{equation}
    m_{\tmop{st}} (\tmmathbf{a} \otimes_{\tmop{st}} \tmmathbf{b} ) =
    m_{\tmop{st}} ( \widetilde{\tmmathbf{a}} \otimes_{\tmop{st}} \tmmathbf{b}
    ) \otimes a_r + m_{\tmop{st}} (\tmmathbf{a} \otimes_{\tmop{st}}
    \widetilde{\tmmathbf{b}} ) \otimes b_s + m_{\tmop{st}} (
    \tilde{\tmmathbf{a}} \otimes_{\tmop{st}} \widetilde{\tmmathbf{b}} )
    \otimes a_r b_s
  \end{equation}
  where $a_r b_s$ is the product in $A$ of $a_r$ and $b_s$.
\end{enumeratenumeric}
For example :
\begin{equation}
  m_{\tmop{st}} ( ( a_1 \otimes a_2 ) \otimes_{\tmop{st}} b_1 ) = a_1 \otimes
  a_2 \otimes b_1 + a_1 \otimes b_1 \otimes a_2 + b_1 \otimes a_1 \otimes a_2
  + a_1 \otimes a_2 b_1
\end{equation}
With this product, $A^{\tmop{st}}$ is a unital algebra (unit $\emptyset$) and
if $A$ is commutative, then $A^{\tmop{st}}$ is commutative. Moreover
\begin{equation}
  \pi_{\tmop{st}} ( A^{\tmop{st}}_{( r )} \otimes_{\tmop{st}} A^{\tmop{st}}_{(
  s )} ) \subset \bigoplus_{t = \max ( r, s )}^{r + s} A^{\tmop{st}}_{( t )}
\end{equation}
On the same way one can define :
\begin{itemizeminus}
  \item a counit $\eta_{\tmop{st}} : A^{\tmop{st}} \rightarrow k$ by
  $\eta_{\tmop{st}} ( \emptyset ) = 1_k$ and for $s \geq 1$, $\eta_{\tmop{st}}
  ( a_1 \otimes \ldots \otimes a_s ) = 0$,
  
  \item a coproduct $\Delta_{\tmop{st}} : A^{\tmop{st}} \rightarrow
  A^{\tmop{st}} \otimes_{\tmop{st}} A^{\tmop{st}}$ such that
  $\Delta_{\tmop{st}} ( \emptyset ) = \emptyset \otimes_{\tmop{st}} \emptyset$
  and for $s \geq 1$ and $\tmmathbf{a}= a_1 \otimes \ldots \otimes a_s \in
  \mathcal{A}^{\tmop{st}}_{( s )}$,
  \begin{equation}
    \Delta_{\tmop{st}} ( \tmmathbf{a} ) = \tmmathbf{a} \otimes_{\tmop{st}}
    \emptyset + \emptyset \otimes_{\tmop{st}} \tmmathbf{a} + \sum_{r = 1}^{s -
    1} ( a_1 \otimes \ldots \otimes a_r ) \otimes_{\tmop{st}} ( a_{r + 1}
    \otimes \ldots \otimes a_s )
  \end{equation}
\end{itemizeminus}
such that $A^{\tmop{st}}$ is a graded coalgebra.

It is a matter of fact to check that $A^{\tmop{st}}$ is a connected bialgebra
(and thus a Hopf algebra) for the filtration :
\begin{equation}
  A^{\tmop{st}}_n = \bigoplus_{k = 0}^n A^{\tmop{st}}_{( n )}
\end{equation}
which is called the stuffle Hopf algebra on $A$ and $A^{\tmop{st}}_0 = k
\emptyset$, ${A^{\tmop{st}}}' = \bigoplus_{n = 1}^{+ \infty} A^{\tmop{st}}_{( n
)}$. We also have, for a sequence $\tmmathbf{a} \in {A^{\tmop{st}}}'$ and $n
\geq 1$,
\begin{equation}
  \Delta'^{[ n ]}_{\tmop{st}} (\tmmathbf{a}) = \sum_{\tmmathbf{a}^1 \otimes
  \ldots \otimes \tmmathbf{a}^n =\tmmathbf{a}} \tmmathbf{a}^1
  \otimes_{\tmop{st}} \ldots \otimes_{\tmop{st}} \tmmathbf{a}^n
\end{equation}
where the sum is over $n$--tuple of non-empty words $(\tmmathbf{a}^1,
\ldots,\tmmathbf{a}^n )$ such that the concatenation of these sequences gives
$\tmmathbf{a}$. In particular, for the antipode, $S ( \emptyset ) = \emptyset$
and if $\tmmathbf{a}= a_1 \otimes \ldots \otimes a_s \in {A^{\tmop{st}}}'$,
\begin{equation}
  \begin{array}{lll}
    S (\tmmathbf{a}) & = & \displaystyle \sum_{k \geq 1} ( - 1 )^k m_{\tmop{st}}^{[ k ]}
    \circ \Delta_{\tmop{st}}'^{[ k ]} (\tmmathbf{a})\\
    & = & \displaystyle\sum_{k \geq 1} ( - 1 )^k \sum_{\tmmathbf{a}^1 \otimes \ldots
    \otimes \tmmathbf{a}^k =\tmmathbf{a}} m_{\tmop{st}}^{[ k ]}
    (\tmmathbf{a}^1 \otimes_{\tmop{st}} \ldots \otimes_{\tmop{st}}
    \tmmathbf{a}^k )
  \end{array}
\end{equation}

If $B$ is an algebra, then once again, there is a convolution on $\mathcal{L}(
A^{\tmop{st}}, B )$ :
\[ \varphi \ast \psi = m_B \circ ( \varphi \otimes \psi ) \circ
   \Delta_{\tmop{st}} \]
and, if $B$ is unital, $(\mathcal{U}( A^{\tmop{st}}, B ), \ast )$ is a group.
Moreover if $B$ is commutative then $\mathcal{C} ( A^{\tmop{st}}, B )$ is a
subgroup.

Finally a map $l \in \mathcal{L}( A, B )$ induces a map $l^{\tmop{st}} \in
\mathcal{U}( A^{\tmop{st}}, B^{\tmop{st}} )$ defined by
\[ l^{\tmop{st}} ( \emptyset ) = \emptyset \hspace{1em} \tmop{and}
   \hspace{1em} l^{\tmop{st}} ( a_1 \otimes \ldots \otimes a_r ) = l ( a_1 )
   \otimes \ldots \otimes l ( a_r ) \hspace{1em} ( r \geq 1 ) \]
and $\Delta_{\tmop{st}} \circ l^{\tmop{st}} = ( l^{\tmop{st}}
\otimes_{\tmop{st}} l^{\tmop{st}} ) \circ \Delta_{\tmop{st}}$. Moreover, if
$l$ is an algebra map, then $l^{\tmop{st}} \circ m_{\tmop{st}} = m_{\tmop{st}}
\circ ( l^{\tmop{st}} \otimes_{\tmop{st}} l^{\tmop{st}} )$, thus
$l^{\tmop{st}}$ is a Hopf morphism.

\subsection{The map $j \in \mathcal{U}( A^{\tmop{st}}, A )$ where $A$ is
unital.}

We shall now illustrate the computations of the previous section on the
following map $j \in \mathcal{C}( A^{\tmop{st}}, A )$ defined by $j (
\emptyset ) = 1_A$, $j ( a_1 ) = a_1$ and $j ( a_1 \otimes \ldots \otimes a_r
) = 0$ if $r \geq 2$. In a sense, this will be the only computation of inverse
and of Birkhoff-Rota-Baxter decomposition we will need.

For the inverse, we get the antialgebra morphism $j^{\ast - 1}$ :
\[ \begin{array}{ccc}
     j^{\ast - 1} & = & u_A \circ \eta_{\tmop{st}} + \sum_{k \geq 1} ( u_A
     \circ \eta_{\tmop{st}} - j )^{\ast^k}
   \end{array} \]
Which means that $j^{\ast - 1} ( \emptyset ) = 1_A$ and for a sequence
$\tmmathbf{a}= a_1 \otimes \ldots \otimes a_s \in {A^{\tmop{st}}}'$,
\begin{equation}
  \begin{array}{ccc}
    j^{\ast - 1} (\tmmathbf{a}) & = &\displaystyle \sum_{k \geq 1} ( - 1 )^k m_A^{[ k ]}
    \circ j^{\otimes_{}^k} \circ {\Delta'}_{\tmop{st}}^{[ k ]} (\tmmathbf{a})\\
    & = &\displaystyle \sum_{k \geq 1} ( - 1 )^k \sum_{\tmmathbf{a}^1 \otimes_{} \ldots
    \otimes \tmmathbf{a}^k =\tmmathbf{a}} m_A^{[ k ]} \circ j^{\otimes^k}
    (\tmmathbf{a}^1 \otimes_{\tmop{st}} \ldots \otimes_{\tmop{st}}
    \tmmathbf{a}^k )\\
    & = &\displaystyle \sum_{k \geq 1} ( - 1 )^k \sum_{\tmmathbf{a}^1 \otimes \ldots
    \otimes \tmmathbf{a}^k =\tmmathbf{a}} j \circ m_{\tmop{st}}^{[ k ]}
    (\tmmathbf{a}^1 \otimes_{\tmop{st}} \ldots \otimes_{\tmop{st}}
    \tmmathbf{a}^k )\\
    & = & ( - 1 )^s a_1 \ldots a_s\\
    & = & j \circ S (\tmmathbf{a})
  \end{array}
\end{equation}

If $( A, p_+ )$ is a Rota-Baxter algebra then the  Bogoliubov preparation map
$\bar{j}$ associated to $j$ is such that $\bar{j} ( \emptyset ) = 1_A$ and is
defined recursively on vector spaces ${A_n^{\tmop{st}}}'$ ($n \geq 1$) by
\begin{equation}
  \bar{j} ( h ) = j ( h ) - m_A \circ ( p_- \otimes \tmop{Id} ) \circ (
  \bar{j} \otimes j ) \circ \Delta_{\tmop{st}}' ( h )
\end{equation}
Let us begin the recursion on the length of the sequence. If $\tmmathbf{a}=
a_1$ then $\bar{j} ( a_1 ) = j ( a_1 ) = a_1$. Now
\begin{equation}
  \bar{j} ( a_1 \otimes a_2 ) = j ( a_1 \otimes a_2 ) - m_A \circ ( p_-
  \otimes \tmop{Id} ) \circ ( \bar{j} \otimes_{} j ) ( ( a_1 )
  \otimes_{\tmop{st}} ( a_2 ) ) = - p_- ( a_1 ) a_2
\end{equation}
and
\begin{equation}
  \begin{array}{ccc}
    \bar{j} ( a_1, a_2, a_3 ) & = & - m_A \circ ( p_- \otimes \tmop{Id} )
    \circ ( \bar{j} \otimes j ) ( ( a_1 \otimes a_2 ) \otimes_{\tmop{st}} (
    a_3 ) )\\
    & = & p_- ( p_- ( a_1 ) a_2 ) a_3
  \end{array}
\end{equation}
Thus, for $r \geq 2$,
\begin{equation}
  \bar{j} ( a_1 \otimes \ldots \otimes a_r ) = - p_- ( \bar{j} ( a_1, \ldots,
  a_{r - 1} ) ) a_r
\end{equation}
It is then easy to prove that

\begin{proposition}
  The Birkhoff-Rota-Baxter decomposition $( j_+, j_- ) \in \mathcal{U}(
  A^{\tmop{st}}, A_+ ) \times \mathcal{U}( A^{\tmop{st}}, A_- )$ such that
  \[ j_- \ast j = j_+ \]
  is given by the formula : for $r \geq 1$ and $\tmmathbf{a}= a_1 \otimes
  \ldots \otimes a_r \in {A^{\tmop{st}}}'$,
  \begin{equation}
    \hspace{1em} \left\{\begin{array}{lllll}
      j_+ (\tmmathbf{a}) & = & p_+ ( \bar{j} (\tmmathbf{a}) ) & = & ( - 1 )^{r
      - 1} p_+ ( p_- ( \ldots ( p_- ( a_1 ) a_2 ) \ldots a_{r - 1} ) a_r )\\
      j_- (\tmmathbf{a}) & = & - p_- ( \bar{j} (\tmmathbf{a}) ) & = & ( - 1
      )^r p_- ( p_- ( \ldots ( p_- ( a_1 ) a_2 ) \ldots a_{r - 1} ) a_r )
    \end{array}\right.
  \end{equation}
  Moreover, if $A$ is commutative then $\mathcal{C}( A^{\tmop{st}}, A_{} )$ is
  a group and $j_{_+}$ and $j_-$ are characters.
\end{proposition}

\begin{proof}
  It remains to prove the last assumption, when $A$ is commutative. Since $j$
  is a character it is sufficient to prove that $j_-$ is a character. By
  induction on $t \in \mathbbm{N}$ we will show that for two sequences
  $\tmmathbf{a}$ and $\tmmathbf{b}$ in $A^{\tmop{st}}$, if $l (\tmmathbf{a}) +
  l (\tmmathbf{b}) = t$, then
  \begin{equation}
    j_- ( m_{\tmop{st}} (\tmmathbf{a} \otimes_{\tmop{st}}
    \text{$\tmmathbf{b}$} ) ) = j_- (\tmmathbf{a}) j_- (\tmmathbf{b})
  \end{equation}
  This identity is trivial for $t = 0$ and $t = 1$ since at least one of the
  sequences is the empty sequence. This also trivial for any $t$ if one of the
  sequence is empty. Now suppose that $t \geq 2$ and that $\tmmathbf{a}= a_1
  \otimes \ldots \otimes a_r \in A^{\tmop{st}}_{( r )}$ and $\tmmathbf{b} =
  b_1 \otimes \ldots \otimes b_s \in A^{\tmop{st}}_{( s )}$ with $r \geq 1$,
  $s \geq 1$ and $r + s = t$. Let $\tilde{\tmmathbf{a}} = a_1 \otimes \ldots
  \otimes a_{r - 1} \in A^{\tmop{st}}_{( r - 1 )}$ ($\tilde{\tmmathbf{a}} =
  \emptyset$ if $r = 1$) and $\tmmathbf{\tilde{b}} = b_1 \otimes \ldots
  \otimes b_{s - 1} \in A^{\tmop{st}}_{( s - 1 )}$ ($\tmmathbf{\tilde{b}} =
  \emptyset$ if $s = 1$), then :
  \[ m_{\tmop{st}} (\tmmathbf{a} \otimes_{\tmop{st}} \tmmathbf{b} ) =
     m_{\tmop{st}} ( \tilde{\tmmathbf{a}} \otimes_{\tmop{st}} \tmmathbf{b} )
     \otimes a_r + m_{\tmop{st}} (\tmmathbf{a} \otimes_{\tmop{st}}
     \tmmathbf{\tilde{b}} ) \otimes b_s + m_{\tmop{st}} ( \tilde{\tmmathbf{a}}
     \otimes_{\tmop{st}} \tmmathbf{\tilde{b}} ) \otimes a_r b_s \]
  Now we have
  \[ j_- (\tmmathbf{a}) = - p_- ( j_- ( \tilde{\tmmathbf{a}} ) a_r ) = - p_- (
     x ) \hspace{1em} \tmop{and} \hspace{1em} j_- (\tmmathbf{b}) = - p_- ( j_-
     ( \tmmathbf{\tilde{b}} ) b_s ) = - p_- ( y ) \]
  Thanks to the Rota-Baxter identity
  \[ \begin{array}{ccc}
       j_- (\tmmathbf{a}) j_- (\tmmathbf{b}) & = & p_- ( x ) p_- ( y )\\
       & = & p_- ( x p_- ( y ) ) + p_- ( p_- ( x ) y ) - p_- ( x y )\\
       & = & p_- ( j_- ( \tilde{\tmmathbf{a}} ) a_r p_- ( j_- (
       \tmmathbf{\tilde{b}} ) b_s ) ) + p_- ( p_- ( j_- ( \tilde{\tmmathbf{a}}
       ) a_r ) j_- ( \tmmathbf{\tilde{b}} ) b_s ) - p_- ( j_- (
       \tilde{\tmmathbf{a}} ) a_r j_- ( \tmmathbf{\tilde{b}} ) b_s )
     \end{array} \]
  but as $A$ is commutative, by induction we get
  \[ \begin{array}{ccc}
       j_- (\tmmathbf{a}) j_- (\tmmathbf{b}) & = & - p_- ( j_- (
       \tilde{\tmmathbf{a}} ) j_- (\tmmathbf{b}) a_r ) - p_- ( j_-
       (\tmmathbf{a}) j_- ( \tmmathbf{\tilde{b}} ) b_s ) - p_- ( j_- (
       \tilde{\tmmathbf{a}} ) j_- ( \tmmathbf{\tilde{b}} ) a_r b_s )\\
       & = & - p_- ( j_- ( m_{\tmop{st}} ( \tilde{\tmmathbf{a}} \otimes_{}
       \tmmathbf{b} ) ) a_r ) - p_- ( j_- ( m_{\tmop{st}} (\tmmathbf{a}
       \otimes \tmmathbf{\tilde{b}} ) ) b_s ) - p_- ( j_- ( m_{\tmop{st}} (
       \tilde{\tmmathbf{a}} \otimes \tmmathbf{\tilde{b}} ) ) a_r b_s )\\
       & = & j_- ( m_{\tmop{st}} ( \tilde{\tmmathbf{a}} \otimes \tmmathbf{b}
       ) \otimes a_r ) + j_- ( m_{\tmop{st}} (\tmmathbf{a} \otimes
       \tmmathbf{\tilde{b}} ) \otimes b_s ) + j_- ( m_{\tmop{st}} (
       \tilde{\tmmathbf{a}} \otimes \tmmathbf{\tilde{b}} ) \otimes a_r b_s )\\
       & = & j_- ( m_{\tmop{st}} ( \tilde{\tmmathbf{a}} \otimes \tmmathbf{b}
       ) \otimes a_r + m_{\tmop{st}} (\tmmathbf{a} \otimes
       \tmmathbf{\tilde{b}} ) \otimes b_s + m_{\tmop{st}} (
       \tilde{\tmmathbf{a}} \otimes \tmmathbf{\tilde{b}} \tmmathbf{} ) \otimes
       a_r b_s )\\
       & = & j_- ( m_{\tmop{st}} (\tmmathbf{a} \otimes \text{$\tmmathbf{b}$}
       ) )
     \end{array} \]
  
\end{proof}

As we will see these formulas are almost sufficient to compute the Birkhoff
decomposition in any connected bialgebra.

\section{The Hopf morphism $\iota : H \rightarrow H'^{\tmop{st}}$.}
\label{sec5}

\begin{theorem}
  Let $H = H_0 \oplus H'$ be a connected bialgebra, then the map $\iota : H
  \rightarrow H'^{\tmop{st}}$ defined by $\iota ( 1_H ) = \emptyset$ and
  \begin{equation}
    \forall h \in H', \hspace{1em} \iota ( h ) = \sum_{k \geq 1} \Delta'^{[ k
    ]} ( h ) \in {{H'}^{\tmop{st}}}'
  \end{equation}
  defines an injective Hopf morphism.
\end{theorem}

\begin{proof}
  This map is well defined since, if $k > n \geq 1$,
  \[ \forall h \in H'_n, \hspace{1em} \Delta'^{[ k ]} ( h ) = 0 \]
  It is obviously linear and injective : For $h_1 = \alpha_1 1_H + p' ( h_1 )$
  and $h_2 = \alpha_2 1_H + p' ( h_2 )$ then, thanks to the graduation of the
  vector space $H'^{\tmop{st}}$, if $\iota ( h_1 ) = i ( h_2 )$ then $\alpha_1
  = \alpha_2$ and $p' ( h_1 ) = p' ( h_2 )$ thus $h_1 = h_2$. This is a
  coalgebra map since
  \[ \Delta_{\tmop{st}} ( \iota ( 1_H ) ) = \Delta_{\tmop{st}} ( \emptyset ) =
     \emptyset \otimes_{\tmop{st}} \emptyset = ( \iota \otimes \iota ) \circ
     \Delta ( 1_H ) \]
  where, to avoid ambiguity, we noted $\otimes_{\tmop{st}}$the tensor product
  of two elements of $H'^{\tmop{st}}$. For $h \in H'$,
  \[ \begin{array}{ccc}
       ( \iota \otimes_{\tmop{st}} \iota ) \circ \Delta ( h ) & = & \displaystyle( \iota
       \otimes_{\tmop{st}} \iota ) ( 1 \otimes h + h \otimes 1 + \Delta' ( h )
       )\\
       & = & \displaystyle\emptyset \otimes_{\tmop{st}} \iota ( h ) + \iota ( h )
       \otimes_{\tmop{st}} \emptyset + \left( \sum_{k \geq 1} \sum_{l \geq 1}
       \Delta'^{[ k ]} \otimes_{\tmop{st}} \Delta'^{[ l ]} \right) \circ
       \Delta' ( h )\\
       & = & \displaystyle\emptyset \otimes_{\tmop{st}} \iota ( h ) + \iota ( h )
       \otimes_{\tmop{st}} \emptyset + \sum_{k \geq 1} \sum_{l \geq 1} (
       \tmop{Id}^{\otimes^k} \otimes_{\tmop{st}} \tmop{Id}^{\otimes^l} ) \circ
       ( \Delta'^{[ k ]} \otimes_{} \Delta'^{[ l ]} ) \circ \Delta' ( h )\\
       & = & \displaystyle\emptyset \otimes_{\tmop{st}} \iota ( h ) + \iota ( h )
       \otimes_{\tmop{st}} \emptyset + \sum_{n \geq 1} \sum_{k = 1}^{n - 1} (
       \tmop{Id}^{\otimes^k} \otimes_{\tmop{st}} \tmop{Id}^{\otimes^l} ) \circ
       \Delta'^{[ n ]} ( h )\\
       & = & \displaystyle\emptyset \otimes_{\tmop{st}} \iota ( h ) + \iota ( h )
       \otimes_{\tmop{st}} \emptyset + \Delta'_{\tmop{st}} ( \iota ( h ) )\\
       & = & \displaystyle\Delta_{\tmop{st}} ( \iota ( h ) )
     \end{array} \]
  But $\iota$ is also an algebra map. Let $g$ and $h$ be two elements of $H$.
  If $g$ or $h$ is in $H_0$ then we get trivially that
  \[ \iota ( g h ) = m_{\tmop{st}} ( \iota ( g ) \otimes_{\tmop{st}} \iota ( h
     ) ) \]
  As in the previous section we will prove by induction on $t \geq 2$ that for
  any positive integer $r$ and $s$ such that $r + s = t$, then
  \[ \forall ( g, h ) \in H'_r \times H'_t, \hspace{1em} \iota ( g h ) =
     m_{\tmop{st}} ( \iota ( g ) \otimes_{\tmop{st}} \iota ( h ) ) \]
  Note that
  \begin{equation}
    \begin{array}{ccc}
      \iota ( h ) & = & \displaystyle\sum_{k \geq 1} \Delta'^{[ k ]} ( h )\\
      & = & \displaystyle h + \sum_{k \geq 1} \Delta'^{[ k + 1 ]} ( h )\\
      & = &\displaystyle h + \sum_{k \geq 1} ( \Delta'^{[ k ]} \otimes \tmop{Id} ) \circ
      \Delta' ( h )\\
      & = &\displaystyle h + \sum_{k \geq 1} ( \Delta'^{[ k ]} \otimes \tmop{Id} ) ( h'_{(
      1 )} \otimes h'_{( 2 )} )\\
      & = & h + \iota ( h'_{( 1 )} ) \otimes h'_{( 2 )}
    \end{array}
  \end{equation}
  \[  \]
  For $t = 2$ ($r = s = 1$) then $\iota ( g ) = g$, $\iota ( h ) = h$ and
  \[ \Delta'^{[ 2 ]} ( g h ) = h \otimes g + g \otimes h \]
  thus
  \[ \iota ( g h ) = g h + h \otimes g + g \otimes h = m_{\tmop{st}} ( ( g )
     \otimes_{\tmop{st}} ( h ) ) \]
  More generally
  \[ \begin{array}{ccc}
       \Delta' ( g h ) & = & h \otimes g + g \otimes h + g h'_{( 1 )} \otimes
       h'_{( 2 )} + h'_{( 1 )} \otimes g h'_{( 2 )}\\
       &  & + g'_{( 1 )} h \otimes g'_{( 2 )} + g'_{( 1 )} \otimes g'_{( 2 )}
       h + g'_{( 1 )} h'_{( 1 )} \otimes g'_{( 2 )} h'_{( 2 )}
     \end{array} \]
  and if $f = g h$
  \[ \begin{array}{ccc}
       \iota ( g h ) & = & f + \iota ( f'_{( 1 )} ) \otimes f'_{( 2 )}
     \end{array} \]
  Now
  \[ \begin{array}{ccc}
       \iota ( g h ) & = & ( g h ) + \iota ( h ) \otimes g + \iota ( g )
       \otimes h + \iota ( g h'_{( 1 )} ) \otimes h'_{( 2 )} + \iota ( h'_{( 1
       )} ) \otimes g h'_{( 2 )}\\
       &  & + \iota ( g'_{( 1 )} h ) \otimes g'_{( 2 )} + \iota ( g'_{( 1 )}
       ) \otimes g'_{( 2 )} h + \iota ( g'_{( 1 )} h'_{( 1 )} ) \otimes g'_{(
       2 )} h'_{( 2 )}\\
       &  & \\
       & = & ( g h ) + h \otimes g + \iota ( h'_{( 1 )} ) \otimes h'_{( 2 )}
       \otimes g + g \otimes h + \iota ( g'_{( 1 )} ) \otimes g'_{( 2 )}
       \otimes h + \iota ( g h'_{( 1 )} ) \otimes h'_{( 2 )}\\
       &  & + \iota ( h'_{( 1 )} ) \otimes g h'_{( 2 )} + \iota ( g'_{( 1 )}
       h ) \otimes g'_{( 2 )} + \iota ( g'_{( 1 )} ) \otimes g'_{( 2 )} h +
       \iota ( g'_{( 1 )} h'_{( 1 )} ) \otimes g'_{( 2 )} h'_{( 2 )}\\
       &  & \\
       & = & ( g h ) + h \otimes g + g \otimes h\\
       &  & + \iota ( h'_{( 1 )} ) \otimes h'_{( 2 )} \otimes g + \iota ( g
       h'_{( 1 )} ) \otimes h'_{( 2 )} + \iota ( h'_{( 1 )} ) \otimes g h'_{(
       2 )}\\
       &  & + \iota ( g'_{( 1 )} ) \otimes g'_{( 2 )} \otimes h + \iota (
       g'_{( 1 )} h ) \otimes g'_{( 2 )} + \iota ( g'_{( 1 )} ) \otimes g'_{(
       2 )} h\\
       &  & + \iota ( g'_{( 1 )} h'_{( 1 )} ) \otimes g'_{( 2 )} h'_{( 2 )}
     \end{array} \]
  By the induction we get
  \[ \begin{array}{ccc}
       \iota ( g h ) & = & m_{\tmop{st}} ( ( g ) \otimes_{\tmop{st}} ( h ) )\\
       &  & + \iota ( h'_{( 1 )} ) \otimes h'_{( 2 )} \otimes g +
       m_{\tmop{st}} ( \iota ( g ) \otimes_{\tmop{st}} \iota ( h'_{( 1 )} ) )
       \otimes h'_{( 2 )} + \iota ( h'_{( 1 )} ) \otimes g h'_{( 2 )}\\
       &  & + \iota ( g'_{( 1 )} ) \otimes g'_{( 2 )} \otimes h +
       m_{\tmop{st}} ( \iota ( g'_{( 1 )} ) \otimes_{\tmop{st}} \iota ( h ) )
       \otimes g'_{( 2 )} + \iota ( g'_{( 1 )} ) \otimes g'_{( 2 )} h\\
       &  & + m_{\tmop{st}} ( \iota ( g'_{( 1 )} ) \otimes_{\tmop{st}} \iota
       ( h'_{( 1 )} ) ) \otimes g'_{( 2 )} h'_{( 2 )}
     \end{array} \]
  As $\iota ( g ) = g + \iota ( g'_{( 1 )} ) \otimes g'_{( 2 )}$ and $\iota (
  h ) = h + \iota ( h'_{( 1 )} ) \otimes h'_{( 2 )}$, we get
  \[ \begin{array}{ccc}
       \iota ( g h ) & = & m_{\tmop{st}} ( ( g ) \otimes_{\tmop{st}} ( h ) )\\
       &  & + \iota ( h'_{( 1 )} ) \otimes h'_{( 2 )} \otimes g +
       m_{\tmop{st}} ( ( g ) \otimes_{\tmop{st}} \iota ( h'_{( 1 )} ) )
       \otimes h'_{( 2 )} + \iota ( h'_{( 1 )} ) \otimes g h'_{( 2 )}\\
       &  & + \iota ( g'_{( 1 )} ) \otimes g'_{( 2 )} \otimes h +
       m_{\tmop{st}} ( \iota ( g'_{( 1 )} ) \otimes_{\tmop{st}} ( h ) )
       \otimes g'_{( 2 )} + \iota ( g'_{( 1 )} ) \otimes g'_{( 2 )} h\\
       &  & + m_{\tmop{st}} ( \iota ( g'_{( 1 )} ) \otimes_{\tmop{st}} \iota
       ( h'_{( 1 )} ) ) \otimes g'_{( 2 )} h'_{( 2 )}\\
       &  & + m_{\tmop{st}} ( ( \iota ( g'_{( 1 )} ) \otimes g'_{( 2 )} )
       \otimes_{\tmop{st}} \iota ( h'_{( 1 )} ) ) \otimes h'_{( 2 )} +
       m_{\tmop{st}} ( \iota ( g'_{( 1 )} ) \otimes_{\tmop{st}} ( \iota (
       h'_{( 1 )} ) \otimes h'_{( 2 )} ) ) \otimes g'_{( 2 )}\\
       & = & m_{\tmop{st}} ( ( g ) \otimes_{\tmop{st}} ( h ) )\\
       &  & + m_{\tmop{st}} ( ( g ) \otimes_{\tmop{st}} ( \iota ( h'_{( 1 )}
       ) \otimes h'_{( 2 )} ) )\\
       &  & + m_{\tmop{st}} ( ( \iota ( g'_{( 1 )} ) \otimes g'_{( 2 )} )
       \otimes_{\tmop{st}} ( h ) )\\
       &  & + m_{\tmop{st}} ( ( \iota ( g'_{( 1 )} ) \otimes g'_{( 2 )} )
       \otimes_{\tmop{st}} ( \iota ( h'_{( 1 )} ) \otimes h'_{( 2 )} ) )\\
       &  & \\
       & = & m_{\tmop{st}} ( \iota ( g ) \otimes_{\tmop{st}} \iota ( h ) )
     \end{array} \]

\end{proof}

This morphisms shows that any connected bialgebra can be canonically
identified to a subalgebra of a stuffle algebra. This will help us to define
$\mathcal{U}( A^{\tmop{st}}, B )$ as a group of function from $\mathcal{U}( H,
A )$ to $\mathcal{U}( H, B )$ where $H$ (resp. $A$, $B$) is a connected
bialgebra (resp. unital algebras).

\section{The map $T :\mathcal{U}( A^{\tmop{st}}, B ) \times \mathcal{U}( H, A
) \rightarrow \mathcal{U}( H, B )$ and associated formulas.} \label{sec6}

\subsection{Definition and properties.}

Let $H$ be a connected bialgebra and $A$, $B$ two unital algebras. For
$\varphi \in \mathcal{U}( H, A )$ and $f \in \mathcal{U}( A^{\tmop{st}}, B )$
we define
\[ T ( f, \varphi ) = f \circ \varphi^{\tmop{st}} \circ \iota \]
It is clear that $T ( f, \varphi )$ is a linear morphism from $H$ to $B$ and
\[ T ( f, \varphi ) ( 1_H ) = f \circ \varphi^{\tmop{st}} \circ \iota ( 1_H )
   = f \circ \varphi^{\tmop{st}} ( \emptyset ) = f ( \emptyset ) = 1_B \]
thus $T ( f, \varphi ) \in \mathcal{U}( H, B )$. Moreover if $\varphi \in
\mathcal{C}( H, A )$ and  $f \in \mathcal{C}( A^{\tmop{st}}, B )$ it is clear,
by composition of algebra morphisms that $T ( f, \varphi ) \in \mathcal{C}( H,
B )$.

There are two fundamental properties :
\begin{enumeratenumeric}
  \item Let $f$ and $g$ in $\mathcal{U}( A^{\tmop{st}}, B )$ and $\varphi$ in
  $\mathcal{U}( H, A )$, then
  \begin{equation}
    \begin{array}{ccc}
      T ( f \ast g, \varphi ) & = & m_B \circ ( f \otimes g ) \circ
      \Delta_{\tmop{st}} \circ \varphi^{\tmop{st}} \circ \iota\\
      & = & m_B \circ ( f \otimes g ) \circ ( \varphi^{\tmop{st}} \otimes
      \varphi^{\tmop{st}} ) \circ \Delta_{\tmop{st}} \circ \iota\\
      & = & m_B \circ ( f \otimes g ) \circ ( \varphi^{\tmop{st}} \otimes
      \varphi^{\tmop{st}} ) \circ ( \iota \otimes \iota ) \circ \Delta\\
      & = & m_B ( T ( f, \varphi ) \otimes T ( g, \varphi ) ) \circ \Delta\\
      T ( f \ast g, \varphi ) & = & T ( f, \varphi ) \ast T ( g, \varphi )
    \end{array}
  \end{equation}
  \item If $A = B$, then $T ( j, \varphi ) = \varphi$ because if $h \in H'$,
  \begin{equation}
    \begin{array}{ccc}
      T ( j, \varphi ) ( h ) & = &\displaystyle j \circ \varphi^{\tmop{st}} \left( h +
      \sum_{k \geq 2} \sum h'_{( 1 )} \otimes \ldots \otimes h'_{( k )}
      \right)\\
      & = &\displaystyle j \left( \varphi ( h ) + \sum_{k \geq 2} \varphi ( h'_{( 1 )} )
      \otimes \ldots \otimes \varphi ( h'_{( k )} ) \right)\\
      & = & \varphi ( h )
    \end{array}
  \end{equation}
\end{enumeratenumeric}

\subsection{The semigroup $(\mathcal{U}( A^{\tmop{st}}, A ), \odot )$.}

We leave the details to the reader but, or $f$ and $g$ in $\mathcal{U}(
A^{\tmop{st}}, A )$ let us define
\[ f \odot g = T ( f, g ) \]
This is a binary operator on $\mathcal{U}( A^{\tmop{st}}, A )$ which is
associative (but non-commutative) and
\[ T ( f, j ) = f \circ j^{\tmop{st}} \circ \iota = f \]
thus $(\mathcal{U}( A^{\tmop{st}}, A ), \odot )$ is semigroup.

\subsection{The inverse in $\mathcal{U}( H, A )$.}

Let $\varphi \in \mathcal{U}( H, A )$. Since $T ( j, \varphi ) = \varphi$, if
$\psi = T ( j^{\ast - 1}, \varphi )$, then
\[ \psi \ast \varphi = T ( j^{\ast - 1}, \varphi ) \ast T ( j, \varphi ) = T (
   j^{\ast - 1} \ast j, \varphi ) = T ( u_A \circ \eta_{\tmop{st}}, \varphi )
   = u_A \circ \eta \]
and the inverse of $\varphi$ is $\psi$. For example, if $h \in H'_3$, then
\[ \iota ( h ) = h + h'_{( 1 )} \otimes h'_{( 2 )} + h'_{( 1 )} \otimes h'_{(
   2 )} \otimes h'_{( 3 )} \]
so,
\[ \varphi^{\tmop{st}} \circ \iota ( h ) = \varphi ( h ) + \varphi ( h'_{( 1
   )} ) \otimes \varphi ( h'_{( 2 )} ) + \varphi ( h'_{( 1 )} ) \otimes
   \varphi ( h'_{( 2 )} ) \otimes \varphi ( h'_{( 3 )} ) \]
and finally
\[ \varphi^{\ast - 1} ( h ) = j^{\ast - 1} \circ \varphi^{\tmop{st}} \circ
   \iota = - \varphi ( h ) + \varphi ( h'_{( 1 )} ) \varphi ( h'_{( 2 )} ) -
   \varphi ( h'_{( 1 )} ) \varphi ( h'_{( 2 )} ) \varphi ( h'_{( 3 )} ) \]
We recover the usual formula for the inverse.

\subsection{The BRB decomposition in $\mathcal{U}( H, A )$.}

Finally, let $\varphi \in \mathcal{U}( H, A )$. Since $T ( j, \varphi ) =
\varphi$, if $\varphi_- = T ( j_-, \varphi )$ and $\varphi_+ = T ( j_+,
\varphi )$, then
\[ \varphi_- \ast \varphi = T ( j_-, \varphi ) \ast T ( j, \varphi ) = T ( j_-
   \ast j, \varphi ) = T ( j_+, \varphi ) = \varphi_+ \]
and, of course, $\varphi_{\pm} \in \mathcal{U}( H, A_{\pm} )$. For example, if
$h \in H'_3$, then
\[ \begin{array}{ccc}
     \varphi_+ ( h ) & = & p_+ ( \varphi ( h ) ) - p_+ ( p_- ( \varphi ( h'_{(
     1 )} ) ) \varphi ( h'_{( 2 )} ) ) + p_+ ( p_- ( p_- ( \varphi ( h'_{( 1
     )} ) ) \varphi ( h'_{( 2 )} ) ) \varphi ( h'_{( 3 )} ) )\\
     \varphi_- ( h ) & = & - p_- ( \varphi ( h ) ) + p_- ( p_- ( \varphi (
     h'_{( 1 )} ) ) \varphi ( h'_{( 2 )} ) ) - p_- ( p_- ( p_- ( \varphi (
     h'_{( 1 )} ) ) \varphi ( h'_{( 2 )} ) ) \varphi ( h'_{( 3 )} ) )
   \end{array} \]

Needless to say that if $A$ is commutative, these computations works in the
subgroup $\mathcal{C}( H, A )$.

\section{Conclusion.}

Once these formulas are given, we get formulas in the different contexts where
renormalization, or rather BRB decomposition is needed :
\begin{itemizeminus}
  \item Renormalization in quantum field theory : the connected Hopf algebra
  is the connected graded Hopf algebra of Feynman graphs. The character is
  given, after dimensional regularization, by Feynman integrals with values
  in a commutative algebra of Laurent series in a parameter $\varepsilon$ : $A
  =\mathcal{A}[ [ \varepsilon ] ] [ \varepsilon^{- 1} ]$ with $A_+
  =\mathcal{A}[ [ \varepsilon ] ]$ and $A_- = \varepsilon^{- 1} \mathcal{A}[
  \varepsilon^{- 1} ]$ (see for example {\cite{ck1}}, {\cite{ck2}},
  {\cite{kur}}).
  
  \item Chen's iterated integrals : this kind  of integrals (including
  multizeta values) define characters on a connected graded Hopf algebras of
  trees or ladders (see for example {\cite{kre}}, {\cite{manchon}}).
  
  \item The Birkhoff decomposition in the group of formal identity-tangent
  diffeomorphism with coefficients in $A =\mathbbm{R}[ [ \varepsilon ] ] [
  \varepsilon^{- 1} ]$. Any element
  \[ f ( x ) = x + \sum_{n \geq 2} f_n ( \varepsilon ) x^n \hspace{1em},
     \hspace{1em} f_n ( \varepsilon ) \in \mathbbm{R}[ [ \varepsilon ] ] [
     \varepsilon^{- 1} ] \]
  can be decomposed : $f_- \circ f = f_+$ with
  \[ \begin{array}{ccccc}
       f_- ( x ) & = & \displaystyle x + \sum_{n \geq 2} f_{-, n} ( \varepsilon ) x^n &  &
       f_{-, n} ( \varepsilon ) \in \varepsilon^{- 1} \mathbbm{R}[
       \varepsilon^{- 1} ]\\
       f_+ ( x ) & = & \displaystyle x + \sum_{n \geq 2} f_{+, n} ( \varepsilon ) x^n &  &
       f_{+, n} ( \varepsilon ) \in \mathbbm{R}[ [ \varepsilon ] ] [
       \varepsilon^{- 1} ]
     \end{array} \]
  This factorization corresponds here to the BRB decomposition in the Faà di
  Bruno Hopf algebra (see {\cite{fig}},{\cite{men}}).

\end{itemizeminus}
The same ideas were also used for the the even-odd factorization of characters
in combinatorial Hopf algebras (see {\cite{ag1}}, {\cite{ag2}} and
{\cite{guo}}).

\end{document}